\documentclass[10pt, a4paper]{amsart}
\usepackage{amsmath}
\usepackage{latexsym}
\usepackage[all]{xy}
\usepackage{color}
\newtheorem{teorema}{Theorem}[section]
\newtheorem{definicion}[teorema]{Definition}
\include{amslatex}
\newtheorem{proposicion}[teorema]{Proposition}
\newtheorem{lema}[teorema]{Lemma}
\newtheorem{corolario}[teorema]{Corollary}

\usepackage[margin=1.4in]{geometry}
\newtheorem{comentario}[teorema]{Remark}

\numberwithin{equation}{section}

\begin{document}
\begin{title}[Quotient rings and fields from the metric geometry point of view]
{Quotient rings and fields of integers from the metric geometry point of view}
\end{title}
\maketitle

\begin{center}
\author{Ricardo Gallego Torrom\'e}\\
\address{}
\bigskip
\end{center}
\begin{abstract}
The theory of Gromov-Hausdorff convergence is applied to sequences of quotient rings of integers, leading to the existence of limit rings as the Gromov-Hausdorff limits. It is also shown that such {\it limit rings} can be endowed with an order relation, although they are not ordered fields. When the construction restricts to quotients $\mathbb{Z}/p\mathbb{Z}$ with $p$ prime, the limits are fields. The relation of those limits with the field of the real numbers $\mathbb{R}$ is discussed, showing that the limit fields are dense in $\mathbb{R}$ but they cannot be identified with $\mathbb{R}$ or with the rational field $\mathbb{Q}$, at least when $\mathbb{R}$ and $\mathbb{Q}$ are endowed with the usual order relations and metric structures, neither they can be identified with the $p$-adic number system. It is also shown that the limit rings and fields have an ordinal larger than the reals and that they are endowed with an Archimedean absolute value function.
\end{abstract}
\bigskip
{\bf MSC Class}: 11H99; 11R21; 51F99; 81Q65.\\
\bigskip
{\bf Key words}: Quotien Rings; Prime Fields; Gromov-Hausdorff Convergence; Emergent Quantum Mechanics.
\bigskip

\section{Introduction}
Gromov-Hausdorff distance provides a powerful theory that allows to compare metric spaces between them as a whole. In particular, the Gromov-Hausdorff convergence finds multiple applications in geometry \cite{BuragoBuragoIvanov,Gromov}, as Gromov's theorem on groups with polynomial growth or on the sudy of limits of Riemannian manifold with non-negative Ricci curvature. More rare are the applications in algebra, with the notable exception of the work on matrix algebras by Rieffel \cite{Rieffel 2004}, motivated by string theory problems.

In this work, quotient rings of integers and their limits in a metric sense are investigated by means of methods of metric geometry. The questions considered are initially motivated in the context of theoretical physics, in particular, from the development of a theory of general dynamical systems for a theory of emergent quantum mechanics, Hamilton-Randers theory \cite{Ricardo general dynamics, Ricardo06,Ricardo2014}. According to such a theory, time parameters are subsets of fields endowed with an order relation, although they are not necessarily ordered fields. Appropriate subsects of the rational numbers $\mathbb{Q}$ or the real numbers $\mathbb{R}$ provide natural candidates for the time parameters of the dynamical systems, but since the theory developed had a general scope, other discrete fields were also considered. Moreover, we were interested in the case when the time parameter sets are endowed with a notion of {\it minimal time lapse} in the sense of a time interval which is minimal and small compared with any macroscopic time interval for all practical purposes and that it can be used in a definition of infinitesimal time differences or in a definition of incremental quotient.

In this context, it is natural to consider the limits $p\to +\infty$ of prime fields  $\mathbb{Z}/p\mathbb{Z}$. If such limits can be established, then they can serve to define time parameter candidates endowed with a natural notion of {\it minimal time lapse} given by the expression $\delta t :=1/(p-1)$ for $p$ prime arbitrarily large. But besides this initial motivation, the existence of such limits of quotient rings and fields can have an intrinsic mathematical interest, since the limits constitute novel examples of rings and fields not yet considered in the literature. Moreover, the quotient rings $\mathbb{Z}/p\mathbb{Z}$ can be endowed with order relations that induce an order relation in the limits rings and fields.

The construction of the limits  $p\to +\infty$ of quotient rings $\mathbb{Z}/p\mathbb{Z}$ or, in the case when $p$  takes values on the set of prime numbers, prime quotient fields, is investigated in this paper. The main technical tool that we apply to show the existence of the limits is the theory of Gromov-Hausdorff convergence as discussed in reference \cite{BuragoBuragoIvanov}, Chapter 7. As a result, we  prove the existence of rings obtained as limits in the sense of Gromov-Hausdorff of sequences of quotient rings endowed with metric structures, while in the case of limits of prime fields, the limits provide examples of non-finite fields. Moreover, such fields are endowed with an order relation but are not ordered fields.

Several questions on the nature and properties of the limit rings and fields that are beyond the metric point of view arise. In particular, the uniqueness of the limits or its dependence in the particular initial sequence is an open question. We show that for the case of successions of prime fields, although the Gromov-Hausdorff limit fields $\mathbb{Z}/\infty \mathbb{Z}$ cannot be identified with $\mathbb{R}$ or with $\mathbb{Q}$, the limit fields are dense in $\mathbb{R}$.

The Gromov-Hausdorff limits of quotient rings and fields admit arbitrarily close approximations in the sense of the Gromov-Hausdorff distance topology by quotients rings $\mathbb{Z}/p\mathbb{Z}$ for large enough $p$. In such a way, the minimal  time lapse $\delta t_p$ can be done arbitrarily small compared with the other scale that appear in the geometry, namely, the diameter of $\mathbb{Z}/p\mathbb{Z}$, that for the metrics used in the theory, it is equal to $1$ for each $p\in \mathbb{N}\setminus\{1\}$. It is also shown that the limit fields are equipped with an Archimedean absolute value map.
Finally, it is illustrated that the models of time parameter based upon our construction of the  limit fields $\mathbb{Z}/\infty \mathbb{Z}$ contain the main properties required for the formulation of differential equations or finite difference equations dynamical models of the class discussed in Hamilton-Randers theory.
\section{Notation and Preliminary considerations}
\subsection*{Definition of the metric in $\mathbb{Z}/p\mathbb{Z}$}
The quotient ring $\mathbb{Z}/p\mathbb{Z}$ is the aggregate of equivalence classes $[k]$ module $p$,
\begin{align*}
[k]_p:=\{n\in\mathbb{Z}\,\,s.t.\,\,n-k=\,q p,\,\,p, \,q\in\,\mathbb{Z}\}
\end{align*}
endowed with the usual sum and multiplication of classes.
For each equivalence class a representant $\,k\in\{0,1,2,...,p-1\}$ can be chosen. For each $p\in\,\mathbb{N}$ there is a natural metric function $d_p:\mathbb{Z}/p\mathbb{Z}\times \mathbb{Z}/p\mathbb{Z}\to \,\mathbb{R}$ defined by the expression
\begin{align}
d_{p}([n]_p,[m]_p):=\,|n_0-m_0|,\,n_0\in\,[n]_p,\,m_0\in\,[m]_p,\,0\leq n_0,m_0\,\leq p-1.
\label{metric function in Zp}
\end{align}
That the function $d_p$ satisfies the axioms of a metric follows from the properties of the usual metric function on $\mathbb{R}$ determined by the modulus function $|,|:\mathbb{R}\to \mathbb{R}^+\cup\{0\}$. The distance associated to the metric function \eqref{metric function in Zp} between $[n]_p$ and $[m]_p$ in $\mathbb{Z}/p\mathbb{Z}$  is  $d_p([n]_p,[m]_p)$ and the absolute value of a class $[n]\in\,\mathbb{Z}/p\mathbb{Z}$ is its distance to the zero class, $\|[n]\|_p:=\,d_p([0]_p,[n]_p)$, determining also the function $\|\cdot \|_p: \mathbb{Z}/p\mathbb{Z}\to \mathbb{R}$.

 The following properties of the distance function $d_p$ are easily established:
 \begin{enumerate}

\item The minimal distance $d_{min}:=\min\{d_p([n]_p,[m]_p),\,[n]_p,[m]_p\in\,\mathbb{Z}/p\mathbb{Z}\}$
is equal to $1$ for any $p\in\,\mathbb{N}$ greater or equal to $2$.

\item For any $p\in\,\mathbb{N} $ greater or equal to $2$, the diameter of $\mathbb{Z}/p\mathbb{Z}$ is the maximum distance between elements in $\mathbb{Z}/p\mathbb{Z}$,
\begin{align*}
diam(\mathbb{Z}/p\mathbb{Z}):=\,\max\{d_p([n]_p,[m]_p),\,[n]_p,[m]_p\in\,\mathbb{Z}/p\mathbb{Z}\}.
\end{align*}
It follows that $diam(\mathbb{Z}/p\mathbb{Z})=\,p-1.$
\end{enumerate}

\subsection*{Normalized metric}
The {\it normalized metric} in $\mathbb{Z}/p\mathbb{Z}$ is defined by the expression
 \begin{align}
 \tilde{d}_p:\mathbb{Z}/p\mathbb{Z} \times \mathbb{Z}/p\mathbb{Z}\to \mathbb{R},\quad ([n]_p,[m]_p)\mapsto \,\frac{1}{diam(\mathbb{Z}/p\mathbb{Z})}\,d_p([n]_p,[m]_p).
 \label{normalized metric}
 \end{align}
 For the normalized metric function $\tilde{d}_p$ the following properties hold good:
 \begin{itemize}
 \item The quotient rings $\mathbb{Z}/p\mathbb{Z}$ have finite diameters,
$ diam_{\tilde{d}_p}(\mathbb{Z}/p\mathbb{Z})=1,\,\,\forall\, p\,\in\,\mathbb{N}\setminus\{1\}$.
 Therefore, the spaces $\{(\mathbb{Z}/p\mathbb{Z},\tilde{d}_p),\,p\in \,\mathbb{N}\}$ are compact but cannot be points (according to the characterization 7.4.6 in \cite{BuragoBuragoIvanov}), since the diameter of each $\mathbb{Z}/p\mathbb{Z}$ is non-zero.

 \item In $\mathbb{Z}/p\mathbb{Z}$ the minimal distance possible by using the metric $\tilde{d}_p$ between elements of $\mathbb{Z}/p\mathbb{Z}$ is given by $1/(p-1)$. Therefore,
 \begin{align*}
 \lim_{p\,\to +\infty}\,\tilde{d}_{min}(\mathbb{Z}/p\mathbb{Z})=\,\lim_{p\to +\infty} \frac{1}{p-1} =\,0.
 \end{align*}

 \end{itemize}
\subsection{Order relation in $\mathbb{Z}/p\mathbb{Z} $} The ring $\mathbb{Z}/p\mathbb{Z}$ can be endowed with an order relation in the following way.
\begin{definicion}
If $[n]_p,\,[m]_p\in\,\mathbb{Z}/p\mathbb{Z}$, we say that
\begin{align*}
[n]_p< \,[m]_p\,\quad \textrm{iff}\,\quad n_0\,< \,m_0,
\end{align*}
$\textrm{with}\, \,n_0\in\,[n]_p,\,m_0\in\,[m]_p,\,0\leq n_0,m_0\leq p-1.$
\label{definicion de orden}
\end{definicion}
The  relation $[n]-p\leq \,[m]_p$ as specified by {\it Definition} \ref{definicion de orden} determines a total order in $\mathbb{Z}/p\mathbb{Z}$.
There are other possible definitions of the same order relation in $\mathbb{Z}/p\mathbb{Z}$, for instance, by taking the analogous definition but based upon other congruent elements in each class. Such alternative definitions lead to the same total order relation than ours.

Let us also remark that the order relation determined by the {\it Definition} \ref{definicion de orden} is not compatible with the arithmetic operations of multiplication and addition of the ring $\mathbb{Z}/p\mathbb{Z}$. However, the order relation defined above is compatible with the notion of absolute value introduced before, because if $d_p([0]_p, [n]_p)<d_p([0]_p,[m]_p)$, then $[n]_p<[m]_p$ and viceversa.
Moreover, one has that
\begin{proposicion}
The absolute value function $\|\cdot \|_p:\mathbb{Z}/p\mathbb{Z}\to \mathbb{R}$ is Archimedean.
\end{proposicion}
\begin{proof}
If the strong triangle inequality
\begin{align}
\|[n]+[m] \|_p\leq\, \max\{\|n \|_p,\|m \|_p\}
\label{strong triangle inequality}
\end{align}
is violated, then the $\|\cdot \|_p$ is Archimedean.
Given the classes $[n]_p =[n_0]_p$ and $[m]_p=[m_0]_p$, with $0\leq n_0,\,m_0\leq p-1$, there are two possible cases:
\begin{align*}
\|[n]+[m] \|_p =\,
\begin{cases}
& n_0+m_0,\,\textrm{if}\, n_0+m_0<p \\
& n_0+m_0-p,\,\textrm{if}\, n_0+m_0 \geq p .
\end{cases}
\end{align*}
It is easy to observe that for each of these two cases the strong triangle inequality is violated.

\end{proof}
 \section{Gromov-Hausdorff limits of quotient rings}
The main result of this work arises from the application of Gromov-Hausdorff convergence theory to the class
\begin{align*}
\mathcal{C}:= \{(\mathbb{Z}/p\mathbb{Z},\tilde{d}_p),\,p\in\, \mathbb{N}\setminus\{1\}\}
\end{align*}
of compact metric spaces with the normalized metrics.
 The Gromov-Hausdorff distance between the metric spaces $(\mathbb{Z}/p_i\mathbb{Z},d_{Z_{p_i}})$ and $(\mathbb{Z}/p_j\mathbb{Z},d_{Z_{p_j}})$ is bounded by applying the notions discussed in \cite{BuragoBuragoIvanov}, section 7.3. Let us consider the functions
 \begin{align}
 \upsilon_p:\mathbb{Z}/p\mathbb{Z}\to \,\mathbb{R},\,[n]_p\mapsto \,n_0,\,s.t.\,0\leq n_0\,\leq\, p-1,\,n_0\in\,[n]_p\subset \,\mathbb{R}.
 \label{isometry of Zp in R}
 \end{align}
 With the standard metric structure on $\mathbb{R}$ provided by the module function $|,|:\mathbb{R}\to\,\mathbb{R}^+\cup\{0\}$, $\upsilon_p$ is an isometric embedding of $\mathbb{Z}/p\mathbb{Z}$ into $\mathbb{R}$. The Hausdorff distance between $\upsilon_p(\mathbb{Z}/p_i\mathbb{Z})\subset\,\mathbb{R}$ and $\upsilon_p(\mathbb{Z}/p_j\mathbb{Z})\subset\,\mathbb{R}$ can be evaluated for each pair of quotients indexed by  $p_i,p_j\,\in {\mathbb{N}}$ and is given by
 \begin{align*}
 d_H(\upsilon_{p_i}(\mathbb{Z}/p_i\mathbb{Z}), \upsilon_{p_j}(\mathbb{Z}/p_j\mathbb{Z}))=\,|p_i-\,p_j|.
  \end{align*}
  From the geometric representation of the above isometric embeddings and from the notion of Gromov-Hausdorff distance, which is upper bounded by a generic Hausdorff distance, the  Gromov-Hausdorff  distance between the compact, metric spaces $(\mathbb{Z}/p_i\mathbb{Z},d_{p_i})$ and $(\mathbb{Z}/p_j\mathbb{Z},d_{p_j})$  is such that
 \begin{align}
 d_{GH}((\mathbb{Z}/p_i\mathbb{Z},d_{Z_{p_i}}), (\mathbb{Z}/p_j\mathbb{Z},d_{Z_{p_j}}))\leq\,|p_i-\,p_j|.
 \end{align}

Let us consider the analogous considerations when the rings $\mathbb{Z}/p\mathbb{Z}$ are endowed with metric structures associated to the normalized distance functions $\tilde{d}_p$ on each $(\mathbb{Z}/p\mathbb{Z},\tilde{d}_p)\in\,\mathcal{C}$. In this case, isometric embeddings are given in the form $\tilde{\upsilon}_p : \mathbb{Z}/p \mathbb{Z}\to \mathbb{R},\,[n]_p\mapsto \,\frac{n_0}{p-1}$, where $n_0$ is the element of $[n]$ in the interval $[0,p-1]$. In this case, the Hausdorff distance is given by the relation
\begin{align*}
 d_H(\tilde{\upsilon}_{p_i}(\mathbb{Z}/p_i\mathbb{Z}), \tilde{\upsilon}_{p_j}(\mathbb{Z}/p_j\mathbb{Z}))=\,\frac{|p_i-p_j|}{\min\{p_i,p_j\}}.
\end{align*}
As a consequence, the Gromov-Hausdorff distance $d_{GH}((\mathbb{Z}/p_i\mathbb{Z},\tilde{d}_{p_i}),(\mathbb{Z}/p_j\mathbb{Z},\tilde{d}_{p_j}))$
associated to $\tilde{d}_p$ is such that the bound
 \begin{align}
 d_{GH}((\mathbb{Z}/p_i\mathbb{Z},\tilde{d}_{p_i}),(\mathbb{Z}/p_j\mathbb{Z},\tilde{d}_{p_j}))\leq\,\frac{|p_i-p_j|}{\min\{p_i,p_j\}}
 \label{Gromov-Hausdorff distance bound}
 \end{align}
 holds good.
The convergence properties associated to the distance functions $\tilde{d}_{GH}$ are more adequate for our purposes than the ones for the distances $d_{GH}$. Therefore, we adopt the distance $\tilde{d}_{GH}$ instead than  ${d}_{GH}$ in the considerations of the theory presented in this paper.

A subset $\mathcal{N}$ of a metric space $X$ is an $\epsilon$-net if $d(x,\mathcal{N})\leq \,\epsilon\,\forall\, x\in X$.
Following Burago et al. (\cite{BuragoBuragoIvanov}, section 7.4) we consider the following notion:
 \begin{definicion}
 A class $\Upsilon$ of compact metric spaces is uniformly totally bounded if the following two conditions hold:
 \begin{enumerate}
 \item There is a constant $\Delta\,> 0$ such that $diam (X)\leq \Delta$ for all $X\in \Upsilon$.

 \item For every $\epsilon\,>0$ there exists a natural number $N(\epsilon)$ such that every $X\in\,\Upsilon$ contains an $\epsilon$-net with no more than $N(\epsilon)$ points.

 \end{enumerate}
 \label{Definition uniformly totally bounded}
 \end{definicion}
 The following result \cite{BuragoBuragoIvanov} is key for our considerations,
 \begin{teorema}
 Any uniformly totally bounded class $\Upsilon$ of compact metric spaces is pre-compact in the Gromov-Hausdorff topology. That is, any sequence of $\Upsilon$ contains a convergent sub-sequence.
 \label{Teorema uniformy totally bounded class}
 \end{teorema}
The uniformly totally bounded property and the corresponding result on convergence of subsequences can
 be applied to sequences of finite rings in $\mathcal{C}$. In particular, let us consider the collection of metric spaces
\begin{align*}
 \mathcal{C}_{\epsilon_0}:=\{(\mathbb{Z}/p\mathbb{Z},\tilde{d}_p)\,s.t.\, p>\,p_0\,\,\textrm{and}\,\,\frac{1}{p-1}<\epsilon_0\}
\end{align*}
Note that for $\epsilon<\,\epsilon_0$ one has that $\mathcal{C}_\epsilon\subset\, \mathcal{C}_{\epsilon_0}$.
Furthermore, we find the following result:
 \begin{lema}
 For every $\epsilon_0> \,0$ the class $\mathcal{C}_{\epsilon_0}$ is uniformly totally bounded.
 \label{tildep is uniformly bounded}
 \end{lema}
 \begin{proof}
 Condition (1) in {\it Definition} \ref{Definition uniformly totally bounded} holds good, since $diam_{\tilde{d}_p}(\mathbb{Z}/p\mathbb{Z})=1$ for all $p\in\,\mathbb{N}\setminus\{1\}$.
 In order to show that also the condition $(2)$ is satisfied, let us consider an arbitrary $\epsilon > 0$ and $1/p<\epsilon\,<\epsilon_0$. Denote by $\mathcal{N}_{p}(\epsilon)$ a regular distribution of $p$ points in $\mathbb{Z}/p'\mathbb{Z}$ and assume an uniform distribution of points on $\mathbb{R}$. In this case, we have for each $[\alpha']\in\,\mathcal{N}_p$ that
     \begin{align*}
     \tilde{d}_{p'}([\alpha]',\mathbb{Z}/p'\mathbb{Z})\leq \,\frac{[p'/p]}{p'-1}=\,\frac{p'-r_0}{p'-1}\,\frac{1}{p},
     \end{align*}
 where $[p'/p]$ is the integer part of the quotient $p'/p$ and $r_0$ depends on $p'$ and $p$ and range in $\{0,...,p-1\}$. It is direct that for a given $\epsilon >0$ we can choose $p$ such that
 \begin{align*}
 \tilde{d}_{p'}([\alpha]',\mathbb{Z}/p'\mathbb{Z})\leq 1/p\leq \epsilon.
 \end{align*}
 This is repeated for every $\mathbb{Z}/p'\mathbb{Z}\in\,\mathcal{C}_{\epsilon_0}$. It follows that
  \begin{align*}
 \tilde{d}_{GH}(\mathcal{N}_p,\mathbb{Z}/p'\mathbb{Z})\leq 1/p\leq \epsilon,\quad \forall \,\mathbb{Z}/p'\mathbb{Z}\in\,\mathcal{C}_{\epsilon_0}.
 \end{align*}
 Thus, given an $0<\epsilon<\epsilon_0$,
from which follows that $\mathcal{N}_{p}(\epsilon)$ is an $\epsilon$-net with $p$ points for each of the elements $\mathbb{Z}/p'\mathbb{Z}\in\,\mathcal{C}_{\epsilon_0}$.
 \end{proof}
It follows the following fundamental result:
 \begin{teorema}
For each $\epsilon_0>0$ the class of compact spaces $\mathcal{C}_{\epsilon_0}$ contains convergent sequences in the Gromov-Hausdorff topology of compact metric spaces.
 \label{fundamental result}
 \end{teorema}
\begin{proof}
  Let us to disregard the first elements of a sequence, whose distance to the limit could be larger than $\epsilon_0$ eliminating from the collection initial elements. For the existence of  a convergent sub-sequence in $\mathcal{C}_{\epsilon_0}$, the initial terms for $p$ small are of no relevance. Then by {\it Lemma} \ref{tildep is uniformly bounded}, {\it Theorem} \ref{Teorema uniformy totally bounded class} is applied to $\mathcal{C}_{\epsilon_0}$.
\end{proof}
\begin{comentario}
This modification of the notion of uniformly bounded family class of metric spaces works in our case because $\epsilon <\epsilon_0$ implies that $\mathcal{C}_\epsilon\subset\, \mathcal{C}_{\epsilon_0}$.
\end{comentario}

 By definition of the Gromov-Hausdorff convergence, if the subsequence $\mathcal{S}=\{\mathbb{Z}/p_n\mathbb{Z},\,n=1,2,3,...\}\subset \,\mathcal{C}$ is convergent, then  the limits
 \begin{align*}
 \lim_{n\to +\infty}\tilde{d}_{p_n}: =\,\tilde{d}_\infty,\quad \lim_{n\to\,+\infty} \,(\mathbb{Z}/p_n\mathbb{Z}) :=\,\mathbb{Z}/\infty\mathbb{Z}
 \end{align*}
exist and the {\it limit metric quotient rings} is
\begin{align*}
(\mathbb{Z}/\infty\mathbb{Z},\tilde{d}_{\infty}):=\,\lim_{n\to +\,\infty}(\mathbb{Z}/p_n\mathbb{Z}, \tilde{d}_{p_n}),
\end{align*}
with $(\mathbb{Z}/p_n\mathbb{Z}, \tilde{d}_{p_n}) \in\,\mathcal{S}\subset \,\mathcal{C}$.
 Let us now consider the map
 \begin{align}
 diam_{\tilde{d}}:\mathcal{C}\to \,\mathbb{R},\,\quad (\mathbb{Z}/p\mathbb{Z},\tilde{d}_{p})\mapsto \,diam_{\tilde{d}_p}(\mathbb{Z}/p\mathbb{Z}).
 \end{align}
  Since the map $diam_{\tilde{p}}:\mathcal{C}\to \mathbb{R}$ is constant, for a given convergent subsequence $\mathcal{S}=\{\mathbb{Z}/p_n\mathbb{Z},\,n=1,2,3,...\}\subset \,\mathcal{C}$ it holds that
\begin{align*}
 diam_{\widetilde{d}_\infty}(\lim_{n\to +\infty}\,\mathbb{Z}/n\mathbb{Z}) =\,\lim_{n\to +\infty} \,diam_{\tilde{d}_{n}}(\mathbb{Z}/n\mathbb{Z})=\,1.
 \end{align*}
Then we have that
 \begin{proposicion}
  For any convergent subsequence $\mathcal{S}\subset\, \mathcal{C}$ convergent to the limit space $(\mathbb{Z}/\infty\mathbb{Z},\tilde{d}_\infty)$, the relation
 \begin{align}
 diam_{\tilde{d}_\infty}\, \mathbb{Z}/\infty\mathbb{Z}=1
 \end{align}
 holds good.
 \label{convergence of diameters}
 \end{proposicion}
Since the diameter  $ diam_{\tilde{d}_\infty}(\mathbb{Z}/\infty\mathbb{Z})$ is finite but it is not zero, we have that
 \begin{corolario}
 Given a convergent subsequence $\mathcal{S}\subset \,\mathcal{C}$, the limit space $\mathbb{Z}/\infty\mathbb{Z}$  is compact but it is not a point.
 \end{corolario}

{\bf Remarks}. First, a direct construction of the set $\mathbb{Z}/\infty\mathbb{Z}$
as an aggregate of limits is not possible by the means considered until now, since the set $\mathbb{Z}/\infty\mathbb{Z}$
contains more elements than any $\mathbb{Z}/p\mathbb{Z}$ of the given sequence $\mathcal{S}$.
Second, one can consider a collection $\tilde{\mathcal{C}}\subset \mathcal{C}$. Then the question of uniqueness of the limits depending on the particular sequence $\tilde{\mathcal{S}}\in\,\tilde{\mathcal{C}}$ arises. In general, one should not expect uniqueness in this sense of the limit processes.
\subsection{Application to sequences of prime fields}
 Let us consider the collection of prime fields as metric spaces,
 \begin{align}
 \widetilde{\mathcal{P}}:=\{(\mathbb{Z}/p_i\mathbb{Z},\tilde{d}_{Z_{p_i}}),\,p_i \in\,\Pi\},
 \end{align}
  where $\Pi\subset \,\mathbb{N}$ denotes the aggregate of the increasing prime numbers.
 Given the collection $\widetilde{\mathcal{P}}$, there exist convergent subsequences in $\widetilde{\mathcal{P}}$ in the Gromov-Hausdorff
 sense,
 \begin{teorema}
The collection $\widetilde{\mathcal{P}}$ contains at least one convergent subsequence of compact metric spaces.
\label{subsequencia primos}
\end{teorema}
The proof mimics the one of {\it Theorem} \ref{fundamental result} and it is not necessary to repeat the argument.

\section{Ring and field structures defined in the limit $\mathbb{Z}/\infty \mathbb{Z}$}
The natural question that arises is related with the existence of algebraic structures in the limit space $\mathbb{Z}/\infty\mathbb{Z}$:
\begin{itemize}
\item Is $\mathbb{Z}/\infty\mathbb{Z} =\lim_{n\to +\infty} S_n$, with $\mathcal{S}\subset\, \mathcal{C}$ a convergent subsequence, furnished with  convenient operations, a ring?

\item Is $\mathbb{Z}/\infty\mathbb{Z}=\lim_{n\to +\infty} S_n$ with $\mathcal{S}\subset\,\widetilde{\mathcal{P}}$ a convergent subsequence, furnished with convenient operations, a field?
\end{itemize}
The following result answer positively the above questions.
\begin{teorema}
Let $(\mathbb{Z}/\infty\mathbb{Z},\tilde{d}_\infty)$ be the limit in the Gromov-Hausdorff sense of a convergent subsequence in $\mathcal{C}\times \mathcal{C}$. Then the limit space $\mathbb{Z}/\infty\mathbb{Z}$ is furnished with a ring structure by the limit operations $+_\infty$ and $\cdot_\infty$ to be defined below. If the convergent subsequence $\mathcal{S}$ is in $\widetilde{\mathcal{P}}$, then $(\mathbb{Z}/\infty\mathbb{Z}, +_\infty, \cdot_\infty)$ is a field.
\label{limit field theorem}
\end{teorema}
\begin{proof}
The algebraic addition operation $+_p:\mathbb{Z}/p\mathbb{Z}\times\,\mathbb{Z}/p\mathbb{Z}\to \mathbb{Z}/p\mathbb{Z}$ is equivalent to the subset of $(\mathbb{Z}/p\mathbb{Z})^3$ of the form
\begin{align*}
\tilde{+}_p:=\{([n]_p,[m]_p,[n+m]_p)\in\,(\mathbb{Z}/p\mathbb{Z})^3,\,[n],[m]\in\,\mathbb{Z}/p\mathbb{Z}\}.
\end{align*}
Let us consider a convergent subsequence $\mathcal{S}=\,\{(\mathbb{Z}/p\mathbb{Z},\tilde{d}_p)\}$ of $\mathcal{C}$ converging to $(\mathbb{Z}/\infty\mathbb{Z},\tilde{d}_\infty)$.
Since $\mathcal{S}$ converges to $(\mathbb{Z}/\infty\mathbb{Z},\tilde{d}_\infty)$ in the Gromov-Hausdorff sense, the sequence $\{(\tilde{+}_p,\bar{d}_p)\}$, with $p$ in the same index set than $\{\mathcal{S}_p\}$ and with distance function defined by
\begin{align*}
\bar{d}_p :=\,\left(\tilde{d}^2_p(1)+\,\tilde{d}^2_p(2)+\tilde{d}^2_p(3)\right)^{1/2}
\end{align*}
 is convergent in the Gromov-Hausdorff sense towards the limit of the sequence $\{\tilde{+}_\infty,\,\bar{d}_\infty)\}$, where
 \begin{align*}
\bar{d}_\infty =\,\left(\tilde{d}^2_\infty(1)+\,\tilde{d}^2_\infty(2)+\tilde{d}^2_\infty(3)\right)^{1/2}
\end{align*}
and where $\tilde{+}_\infty\,\subset \,\left(\mathbb{Z}/\infty\mathbb{Z}\right)^3$.
The limit $\{(\tilde{+}_\infty,\,\bar{d}_\infty)\}$ determines the addition  $+_\infty:\mathbb{Z}/\infty\mathbb{Z}\times \,\mathbb{Z}/\infty\mathbb{Z}\to \mathbb{Z}/\infty\mathbb{Z}$ by the action of the canonical projections
\begin{align*}
pr_i:\left(\mathbb{Z}/\infty \mathbb{Z}\right)^3\to \mathbb{Z}/\infty\mathbb{Z}(i),\quad i=1,2,3
\end{align*}
when restricted to the subset $\tilde{+}_\infty\subset \,\left(\mathbb{Z}/\infty\mathbb{Z}\right)^3$,
\begin{align*}
+_\infty:& pr_1(\tilde{+}_\infty)\times \,pr_2(\tilde{+}_\infty)\to pr_3(\tilde{+}_\infty),\quad (pr_1(\alpha)\times \,pr_2(\alpha))\mapsto pr_3(\alpha).
\end{align*}
Since $ pr_i(\tilde{+}_\infty)\simeq \,\mathbb{Z}/\infty\mathbb{Z}$, the aggregate $\tilde{+}_\infty$ and the projections $pr_i$ determine an abelian group structure in $\mathbb{Z}/\infty\mathbb{Z}$. To show this, one starts noting that since $[n]_p+\,[m]_p =\,[m]_p+\,[n]_p$ for all $p\in \mathbb{N}\setminus\{1\}$, then the commutativity also hold for the limit operation. Otherwise there will be two limits for the same convergent subsequence. The rest of the group axioms for $(\mathbb{Z}/p\mathbb{Z},+_{\infty})$ can be proved by analogous arguments.

Similarly, the product operation $\cdot_\infty:\mathbb{Z}/\infty\mathbb{Z}\times \mathbb{Z}/\infty\mathbb{Z}\to \mathbb{Z}/\infty\mathbb{Z}$ is constructed using the the limits of subsets of the product
\begin{align*}
\tilde{\cdot}_p:=\{([n]_p,[m]_p,[n\cdot m]_p)\in\,(\mathbb{Z}/p\mathbb{Z})^3,\,[n]_p,[m]_p\in\,\mathbb{Z}/p\mathbb{Z}\}.
\end{align*}
By taking limits in a similar way as before, one obtains $\tilde{\cdot}_\infty=\,\lim_{p\to +\infty}\tilde{\cdot}_p\,\subset \, \left(\mathbb{Z}/\infty\mathbb{Z}\right)^3$. Then the projections $pr_i$ and the corresponding limits $\tilde{\cdot}_\infty$ determine a product operation on $\mathbb{Z}/\infty\mathbb{Z}$,
\begin{align*}
\tilde{\cdot}_{\infty}:pr_1(\tilde{\cdot}_\infty)\times \,pr_2(\tilde{\cdot}_\infty)\to pr_3(\tilde{\cdot}_\infty),\quad (pr_1(\alpha)\times \,pr_2(\alpha))\mapsto pr_3(\alpha).
\end{align*}

Finally, in the case of convergent sequences in $\widetilde{\mathcal{P}}$, the limit inverse operation $^{-1}_\infty$ is also well defined (except for the zero element), endowing the limit $\mathbb{Z}/\infty \mathbb{Z}$ with a field structure.
\end{proof}

\section{$\epsilon$-Approximations}
We introduce the notion of $\epsilon$-approximation, a notion useful to prove the existence of an order relation in $\mathbb{Z}/\infty\mathbb{Z}$.
 \begin{definicion}
Let $\mathcal{S}$ be a subsequence of $\mathcal{C}_{\epsilon_0}$ convergent to $(\mathbb{Z}/\infty\mathbb{Z},\tilde{d}_\infty)$. Then $(\mathbb{Z}/p\mathbb{Z},\tilde{d}_p)\in\,\mathcal{S}$ is an $\epsilon$-approximation to $(\mathbb{Z}/\infty\mathbb{Z},\tilde{d}_\infty)$ if
$d_{GH}((\mathbb{Z}/p\mathbb{Z},\tilde{d}_p),(\mathbb{Z}/\infty\mathbb{Z},\tilde{d}_\infty))<\epsilon.$
\label{epsilon approximation}
 \end{definicion}
 This notion of $\epsilon$-approximation is essentially the same as the one in \cite{BuragoBuragoIvanov}, as one can deduce from {\it Theorem 7.4.11} there.

It is direct that if the sequence converges $\mathcal{S}\to \,(\mathbb{Z}/\infty\mathbb{Z},\tilde{d}_\infty)$, then for every $\epsilon >0$ there is a $p\in \mathbb{N}$ such that $(\mathbb{Z}/p\mathbb{Z}, \tilde{d}_p)\in \mathcal{S}$ is an $\epsilon$-approximation.

Let us consider now the order relations
\begin{align*}
[n]_{p_1}\leq \,[m]_{p_1},\,
[n]_{p_2}\leq \,[m]_{p_2},\,[n]_{p_3}\leq \,[m]_{p_3},\, ...
\end{align*}
where $n,m<\,p_1<\,p_2<\,p_3\,<...$ are such that $(\mathbb{Z}/p_i\mathbb{Z},\tilde{d}_{p_i})\in\,\mathcal{S}$ and $[n]_{p_i},\,[m]_{p_i}\in\,\mathbb{Z}/p_i\mathbb{Z}$. For every $\epsilon >0$ there is an $\epsilon$-approximation of $(\mathbb{Z}/\infty\mathbb{Z},\tilde{d}_\infty)$ such that the  relation $[n]_p\leq \,[m]_p$ is valid for any $p>\,p_0(\epsilon)$, for a large enough $p_0$. This condition  is equivalent to define a function which is constant for $p>p_0(\epsilon)$ and hence, continuous for large $p$, that allows to define an order relation in the limit space $(\mathbb{Z}/\infty\mathbb{Z}, \tilde{d}_\infty)$ as the value of the limits of all possible characteristic functions $f(n,m)_p$ for $p\to+\infty$ and all possible $n,m\in\mathbb{N}\setminus\{1\}$, where the function $f([m],[n])_p$ is defined by
\begin{align*}
f([m],[n])_p:=&\,1,\,\textrm{if}\, [n]_p\geq [m]_p,\,[n_p],[m]_p\in\mathbb{Z}/p\mathbb{Z},\\
&\,0,\,\textrm{if}\, [m]_p> [n]_p,\,[n]_p,[m]_p\in\mathbb{Z}/p\mathbb{Z}.
\end{align*}
The sequence of functions $\{f_p\}$ is convergent in the limit $p\to \,+\infty$.
As a consequence of this construction, we have that $f:\mathbb{Z}/\infty\mathbb{Z}\times \,\mathbb{Z}/\infty\mathbb{Z}\to\,\mathbb{R}$ determines a total order relation in the limit ring $\mathbb{Z}/\infty\mathbb{Z}$:
\begin{definicion}
Two classes $[m]_\infty,[n]_\infty\in\,\mathbb{Z}/\infty \mathbb{Z}$ are related by $[m]_\infty\leq\,[n]_\infty$  there is an $\epsilon$-approximation such that $[m]_p\leq \,[n]_p$ for any $p$ large enough in the sub-sequence convergent $\{S_p\}$.
\end{definicion}

The above arguments imply the following result,
\begin{teorema}
The limit ring  $\mathbb{Z}/\infty\mathbb{Z}$ can be endowed with a total order.
\label{ordered field}
\end{teorema}
Let us remark that although the limit field $\mathbb{Z}/\infty\mathbb{Z}$ admits a total ordered relation, it is not an ordered ring, since the algebraic operations are not compatible with the order relation.
\begin{proposicion}
The absolute value function
 $\|\cdot\|_\infty:\mathbb{Z}/\infty\mathbb{Z}\times \,\mathbb{Z}/\infty\mathbb{Z}\to \,\mathbb{R}$ is Archimedean .
\end{proposicion}
\begin{proof}
Since each $\epsilon$-approximation  $\|\cdot \|_p :\mathbb{Z}/p\mathbb{Z}\to \mathbb{R}$ is Archimedean, then the same holds for the limit  $\|\cdot \|_\infty :\mathbb{Z}/p\mathbb{Z}\to \mathbb{R}$
\end{proof}
Thus the system $(\mathbb{Z}/\infty\mathbb{Z}$ cannot be identify with the $p$-adic number system with the $p$-adic norm.

The notion of $\epsilon$-approximation also allows to introduce the notion of infinitesimal distance in $\mathbb{Z}/\infty\mathbb{Z}$. Two elements $[n]_\infty,[m]_\infty\in \,\mathbb{Z}/\infty\mathbb{Z}$ are infinitesimally close if for every $\epsilon >0$ there is an $\epsilon$-approximation such that $\tilde{d}_{p}([n]_p,[m]_p)=\,1/(p-1)<\,\epsilon$. The infinitesimal lapse of time is determined by the set of values $\{1/(p-1)\}$ such that are less than a required $\epsilon$.
\section{The relation of the limit fields with the real field $\mathbb{R}$}
\begin{proposicion}
Each of the limits $\mathbb{Z}/\infty\mathbb{Z}$ is dense in the reals $\mathbb{R}$.
\end{proposicion}
\begin{proof}
 Let us consider the injective map
$ \varphi_p:\mathbb{Z}/p\mathbb{Z}\to \mathbb{S}^1\subset \mathbb{C},\,[n]_p\mapsto \,\exp (2\pi\,\imath n/p)$.
 Let $\mathcal{S}$ be a convergent sequence of metric spaces $(\mathbb{Z}/p\mathbb{Z},\tilde{d}_p)$ in the Gromov-Hausdorff sense.
 The collection of maps $\{\varphi_p,\,s.t.\,(\mathbb{Z}/p\mathbb{Z},\tilde{d}_p)\in\,\mathbb{S}\}$ determines the limit map $ \varphi_\infty:\mathbb{Z}/\infty\mathbb{Z}\to \mathbb{S}^1$.  Setting $p$ large enough, for each point  $z\in\,\mathbb{S}^1 $ there is a point of $[n]_p(z)\in\,\mathbb{Z}/p\mathbb{Z}$ whose image $\varphi_p([n]_p(z))$ is close to $z$ by a distance less than a given $\delta> \,0$. Then by applying the inverse stereographic projection $s^{-1}:\mathbb{S}^1\setminus\{(0,1)\}\to \mathbb{R}$ to the points $s^{-1}(\varphi_p ([n]_p(z)))$ and $s^{-1}(z)$, by a suitable choice of $\delta$, the images became arbitrarily close to each other. Since $s^{-1}$ is an homeomorphism, then for every $y\in \mathbb{R}$ there is a $p_0\in\,\mathbb{N}$ such that for every $p\geq \,p_0$ there is a point $[n]_p\in\,\mathbb{Z}/p\mathbb{Z}$ whose image by $s^{-1}\circ \varphi_p$ is arbitrarily close to $y\in\,\mathbb{R}$ with the standard topology of $\mathbb{R}$. By continuity, this conclusion is generalized to the limit $s^{-1}\circ \varphi_\infty ([n]_\infty)$ and to any point $y\in\,\mathbb{R}$.
 \end{proof}

Despite this result, we have that
 \begin{proposicion}
 The ordinal of $\mathbb{Z}/\infty\mathbb{Z}$ is larger than the ordinal of the real system.
 \end{proposicion}
 \begin{proof}
 Let us consider the extension of the real numbers by attaching a point representing the infinite extremes of the real line, $\mathbb{R}\cup\infty$ and also its compatification via the extension by the stereographic projection $s:\mathbb{R}\cup \{\infty \}\to \mathbb{S}^1$ such that $s((0,1))=\infty$. Thus, $\mathbb{S}^1$ is a bijective compactification of the extended real line $\mathbb{R}\cup\{\infty\}$. On the other hand, the limit when $p\to +\infty$ of the images $\varphi_p (\mathbb{Z}/p\mathbb{Z})$ is $\mathbb{S}^1$, but there is not injective extension  of  $\varphi_{\infty}:\mathbb{Z}/\infty\mathbb{Z}\to \,\mathcal{S}^1$  to $\mathbb{R}\cup\{\infty\}$, since $(0,1)\in \mathbb{S}$ is the image of $[0]$, while $(0,0)\in \,\mathbb{S}^1$ is the image of $\left[int\{\frac{p+1}{2}\}\right]$ for all $p$, where $int\{\frac{p+1}{2}\}$ is the integer part; all other elements of $\mathbb{S}^1$ are image of elements of $\mathbb{Z}/\infty\mathbb{Z}$ by $\varphi_\infty$. Since $\varphi_\infty$ and $s$ are bijective maps, then $\mathbb{R}\cup \{\infty\}$ and $\mathbb{Z}/\infty\mathbb{Z}$ have the same order.
 \end{proof}
  As a direct consequence we have that,
   \begin{corolario}
   the metric space $(\mathbb{Z}/\infty\mathbb{Z},\tilde{d}_\infty)\neq \,(\mathbb{R},|\,|)$ with the standard metric distance determined by the modulus function $|\,|:\mathbb{R}\to \mathbb{R}$.
   \end{corolario}
    This conclusion is also reached from the point of views of metric spaces, since $(\mathbb{Z}/\infty\mathbb{Z},\tilde{d}_\infty)$ is compact while $(\mathbb{R},||)$ is non-compact. Therefore, for the limits of quotients prime fields, the limits $(\mathbb{Z}/\infty\mathbb{Z},\tilde{d}_\infty)$ constitute a new class of fields endowed with metric functions and order relations, dense in $\mathbb{R}$ and where the order relation is not equivalent to the order relation of the real number system. Hence $\mathbb{Z}/\infty\mathbb{Z}$ cannot be identified  with $\mathbb{R}$ or with the rational number system $\mathbb{Q}$.

 \section{Application to Emergent Quantum Mechanics}
The existence of two distance scales in the fields $\mathbb{Z}/\infty \mathbb{Z}$, $1/p$ and $p$ leads to an interpretation compatible with the two-time dynamics description of quantum systems that appears in Hamilton-Randers theory. When the quotient ring represents a fundamental cycle of Hamilton-Randers theory, $1/p$ is the scale of the minimal internal or microscopic time lapse, being the {\it internal time parameter} identified with each of the classes $[n]_p\in\,\mathbb{Z}/p\mathbb{Z}$, while $1$ is the scale of the minimal macroscopic time lapse. This interpretation links the elements of $\mathbb{Z}/\infty \mathbb{Z}$ with the states of an infinite cellular automaton system, while the analogous link is drawn for the $\epsilon$-approximations for finite $p$. In this setting, an evolution law is a map $U:\mathbb{Z}/\infty \mathbb{Z}\times \,\mathbb{Z}/\infty \mathbb{Z}\to \mathbb{Z}/\infty \mathbb{Z}$ such that it is an homomorphism with respect the product in the first entry. 

In the limit $p\to \infty$, the $\epsilon$-approximation construction offers a method to provide arbitrarily small internal time lapses, leading to an effective use of continuous models for the fundamental dynamics that uses the internal time parameter. Moreover, if $p$ takes prime values, in order to assure the possibility to solve algebraic equations associated to the dynamic, relates the theory of Hamilton-Randers dynamical systems with the prime numbers: the larger is the physical system in question, the larger is the $p$ of the $\epsilon$-approximation. As a grosse description, $p$ large corresponds to large mass systems used in the definition of the time parameter. $\mathbb{Z}/p\mathbb{Z}$ is related with the description of a physical system that can be used to define a clock (time parameter) for the fundamental dynamics as assumed in Hamilton-Randers theory in such a way that the larger the prime $p$, the larger the associated inertial mass \cite{Ricardo2014, Ricardo2023}.

\section{Conclusion and Discussion}
It has been shown that Gromov-Hausdorff convergence leads to the construction of the limits of quotient rings and fields. The key point to apply metric geometry to algebra as we did rest on the consideration of the metric functions $d_p:\mathbb{Z}/p\mathbb{Z}\times \mathbb{Z}/p\mathbb{Z}\to \,\mathbb{R}$ given by \eqref{metric function in Zp}, at the expenses of losing the condition of ordered compatible with the ring operations. Such mathematical structures allow the possibility to apply the theory of Gromov-Hausdorff convergence of metric spaces, with the result of the existence of limit quotient rings and fields.
The limit fields $\mathbb{Z}/\infty \mathbb{Z}$ cannot be identified with $\mathbb{R}$ or $\mathbb{Q}$ because they are ordered fields, while the limit fields are not. That is, despite that the fields $\mathbb{Z}/\infty \mathbb{Z}$ are endowed with a total order relation compatible with the distance function, the order relation is not compatible with the field algebraic operations. Moreover, the ordinal of $\mathbb{Z}/\infty \mathbb{Z}$ is strictly larger than the ordinal of $\mathbb{R}$.  Neither it is possible to identify $\mathbb{Z}/\infty\mathbb{Z}$ with the $p$-adic system of numbers, since this last system is non-Archimedean.

 In connection with its potential use of the limits $\mathbb{Z}/\infty\mathbb{Z}$ as models of time in Hamilton-Randers theory. the notion of $\epsilon$-approximation, that provides an operative notion of minimal lapse of time in such a context of quotient rings and fields, appears rather natural, since it allows a notion of {\it infinitesimal time lapse}, reinforcing the status for time parameter spaces for the dynamical systems used in emergent quantum mechanics \cite{ Ricardo06, Ricardo2014} as alternative to the conventional time parameters as subsets of the reals \cite{Ricardo general dynamics}.

The identification of the elements of $\mathbb{Z}/\infty\mathbb{Z}$ with the values of the time parameters leads to a notion of $2$-dimensional time, reinforced by the appearance of two scales: given an $\epsilon$-approximation, $diam (\mathbb{Z}/p\mathbb{Z} )=1$ provides a scale of time (macroscopic time), while $1/p$ is the scale for the internal time. According to such interpretation,, the notion of time parameters based upon the limit fields $\mathbb{Z}/\infty\mathbb{Z}$ and the $\epsilon$-approximations defines a model for a two-dimensional time.

It appears as unnatural that the construction of a time parameter for a fundamental dynamics should rely on high level analytical and geometric methods, that rely themselves on the existence of the real number field. Hence the problem of understanding the nature of the limit fields or a construction method of them that do not requires the use of the real number system appears as an urgent question from the point of view of its relevance for modeling fundamental dynamical systems. A second question that remains open in our approach is on the uniqueness of the limits fields and rings, since in principle they can depend upon the particular initial sequence $\tilde{\mathcal{C}}$.

\end{document}